\documentclass[10pt]{amsart}

\usepackage{amsthm,amsfonts,amsmath,verbatim,amssymb,verbatim,cite,hyperref}

\numberwithin{equation}{section}
\newtheorem{theorem}{Theorem}[section]
\newcommand{\qbin}[2]{\genfrac{[}{]}{0pt}{}{#1}{#2}}

\newcommand{\abs}[1]{\lvert#1\rvert}
\newcommand{\la}{\lambda}
\DeclareMathOperator{\vs}{vs}
\DeclareMathOperator{\hs}{hs}
\DeclareMathOperator{\sk}{sk}
\DeclareMathOperator{\ks}{ks}

\begin{document}

\title[Skew Pieri rules]{Remarks on the paper ``Skew Pieri rules for 
Hall--Littlewood functions'' by Konvalinka and Lauve}

\author{S. Ole Warnaar}
\address{School of Mathematics and Physics,
The University of Queensland, Brisbane, QLD 4072, Australia}

\thanks{Work supported by the Australian Research Council}


\begin{abstract}
In a recent paper Konvalinka and Lauve proved several skew Pieri rules 
for Hall--Littlewood polynomials. In this note we show that $q$-analogues 
of these rules are encoded in a $q$-binomial theorem 
for Macdonald polynomials due to Lascoux and the author.
\end{abstract}

\maketitle

\section{The Konvalinka--Lauve formulas and their $q$-analogues}
We refer the reader to \cite{Macdonald95} for definitions 
concerning Hall--Littlewood and Macdonald polynomials.

Let $P_{\la/\mu}=P_{\la/\mu}(X;t)$ and $Q_{\la/\mu}=Q_{\la/\mu}(X;t)$ 
be the skew Hall--Littlewood polynomials, $e_r=P_{(1^r)}$ the $r$th 
elementary symmetric function, $h_r$ the $r$th complete symmetric 
function and $q_r=Q_{(r)}$.
Then the ordinary Pieri formulas for Hall--Littlewood polynomials 
are given by \cite{Macdonald95}
\begin{subequations}
\begin{align} \label{P0a}
P_{\mu} e_r &=\sum_{\la} \vs_{\la/\mu}(t) P_{\la} \\
P_{\mu} q_r &=\sum_{\la} \hs_{\la/\mu}(t) P_{\la} ,
\label{P0b}
\end{align}
\end{subequations}
where the sums on the right are over partitions $\la$ such that
$\abs{\la}=\abs{\mu}+r$.
The Pieri coefficient $\vs_{\la/\mu}(t)$ is given by 
\cite[p. 215, (3.2)]{Macdonald95}
\begin{equation}\label{vs}
\vs_{\la/\mu}(t)=\prod_{i\geq 1} \qbin{\la'_i-\la'_{i+1}}{\la'_i-\mu'_i}_t,
\end{equation}
so that $\vs_{\la/\mu}(t)$ is zero unless $\mu\subseteq\la$ with $\la-\mu$
a vertical $r$-strip.
Similarly, $\hs_{\la/\mu}(t)$ vanishes unless $\mu\subseteq\la$ with
$\la-\mu$ a horizontal $r$-strip, in which case
\cite[p. 218, (3.10)]{Macdonald95}
\begin{equation}\label{hs}
\hs_{\la/\mu}(t)=\prod_{\substack{\la'_i=\mu'_i+1 \\ \la'_{i+1}=\mu'_{i+1}}} 
\big(1-t^{\la'_i-\la'_{i+1}}\big).
\end{equation}

To express the skew Pieri formulas, Konvalinka and Lauve \cite{KL12}
(see also \cite{Konvalinka12}) introduced a third Pieri coefficient
\begin{equation}\label{sk}
\sk_{\la/\mu}(t):=t^{n(\la/\mu)} \prod_{i\geq 1} 
\qbin{\la'_i-\mu'_{i+1}}{\la'_i-\mu'_i}_t,
\end{equation}
where $n(\la/\mu):=\sum_{i\geq 1} \binom{\la'_i-\mu'_i}{2}$.
Note that $\sk_{\la/\mu}(t)=0$ if $\mu\not\subseteq\la$.

It seems Konvalinka and Lauve have been unaware that the above function
has appeared in the literature before.
Indeed, the right-hand side is exactly \cite[Equation (4.3)]{WZ12}, which is
a formula for the modified Hall--Littlewood polynomial
$Q'_{\la/\mu}(1)=Q_{\la/\mu}(1,t,t^2,\dots)$---a result which
first appeared in \cite[Theorem 3.1]{Lascoux05}, albeit in the 
not-so-easily-recognisable form
\[
Q'_{\la/\mu}(1)=
\begin{cases} \displaystyle 
t^{n(\la/\mu)}
\prod_{i=1}^{l(\mu)}
\frac{1-t^{\la'_{\mu_i-i+1}}}{(t;t)_{\mu'_i-\mu'_{i+1}}}
& \text{for $\mu\subseteq\la$}, \\[4mm]
0 & \text{otherwise},
\end{cases}
\]
and in the two papers \cite[p. 173, Remark 2]{Rains06} and 
\cite[Proposition 3.2]{W05} in a more general form pertaining to
Macdonald polynomials, see \eqref{RW} below.
Prior to these four papers the expression for $\sk_{\la/\mu}(t)$
appeared in the theory of abelian $p$-groups:
\[
\sk_{\la/\mu}(t)=t^{n(\la)-n(\mu)}\alpha_{\la}(\mu;t^{-1}),
\]
where $\alpha_{\la}(\mu;p)$ is the number of subgroups of type $\mu$
in a finite abelian $p$-group of type $\la$,
\cite{Butler94,Delsarte48,Dyubyuk48,Yeh48}.

\begin{theorem}[{Konvalinka--Lauve \cite[Theorems 2--4]{KL12}}]\label{thmKL}
For partitions $\nu\subseteq\mu$,
\begin{subequations}
\begin{align}\label{P1}
P_{\mu/\nu} e_r&=\sum_{\la,\eta}(-1)^{\abs{\nu-\eta}}
\vs_{\la/\mu}(t) \sk_{\nu/\eta}(t) P_{\la/\eta} \\
P_{\mu/\nu} h_r&=\sum_{\la,\eta}(-1)^{\abs{\nu-\eta}}
\sk_{\la/\mu}(t) \vs_{\nu/\eta}(t) P_{\la/\eta}
\label{P2} \\
P_{\mu/\nu} q_r&=\sum_{\la,\eta,\omega} 
(-1)^{\abs{\nu-\omega}} t^{\abs{\omega-\eta}}  
\hs_{\la/\mu}(t) \vs_{\nu/\omega}(t)\sk_{\omega/\eta}(t) P_{\la/\eta},
\label{P4}
\end{align}
\end{subequations}
where each of the multiple sums is subject to the restriction
$\abs{\la}+\abs{\eta}=\abs{\mu}+\abs{\nu}+r$.
\end{theorem}
For $\nu=0$ the first and third skew Pieri formulas reduce to \eqref{P0a}
and \eqref{P0b} respectively, whereas the second formula simplifies 
to \cite[Theorem 1]{KL12}
\[
P_{\mu} h_r=\sum_{\la} \sk_{\la/\mu}(t)P_{\la/\eta}.
\]
Theorem~\ref{thmKL} for $t=0$ gives the skew Pieri rules for 
Schur functions due to Assaf and McNamara \cite{AM11}
who, more generally, conjectured a skew Littlewood--Richardson rule.
The identities \eqref{P1} and \eqref{P2} were first conjectured by
Konvalinka in \cite{Konvalinka12}. The subsequent proof of
the theorem by Konvalinka and Lauve combines 
Hopf algebraic techniques in the spirit of the proof
of the Assaf--McNamara conjecture \cite{LLS11}
with intricate manipulations involving $t$-binomial coefficients.

The aim of this note is to point out that all of the skew Pieri 
formulas \eqref{P1}--\eqref{P4} are implied by a generalised 
$q$-binomial theorem for Macdonald polynomials and, 
consequently, have simple $q$-analogues.

From here on let $P_{\la/\mu}=P_{\la/\mu}(X;q,t)$ and 
$Q_{\la/\mu}=Q_{\la/\mu}(X;q,t)$ denote skew Macdonald polynomials.
Let $f$ be an arbitrary symmetric function.
Adopting plethystic or $\la$-ring notation, see e.g., 
\cite{Haglund08,Lascoux03}, we define
$f\big((a-b)/(1-t)\big)$
in terms of the power sums with positive index $r$ as
\[
p_r\Big(\frac{a-b}{1-t}\Big)=\frac{a^r-b^r}{1-t^r}.
\]
In other words, $p_r\big((a-b)/(1-t)\big)=a^r \epsilon_{b/a,t}(p_r)$
with $\epsilon_{u,r}$ Macdonald's evaluation homomorphism
\cite[p. 338, (6.16)]{Macdonald95}.
Equivalently, in terms of complete symmetric functions,
\[
h_r\Big(\frac{a-b}{1-t}\Big)=[z^r]\,\frac{(bz;t)_{\infty}}{(az;t)_{\infty}}.
\]
We now define the following five Pieri coefficients for Macdonald 
polynomials:
\begin{subequations}\label{Pcoeff}
\begin{align}
\vs_{\la/\mu}(q,t)&:=\psi'_{\la/\mu}(q,t)=
(-1)^{\abs{\la-\mu}} Q_{\la/\mu}\Big(\frac{q-1}{1-t}\Big) \\
\hs_{\la/\mu}(q;t)&:=\varphi_{\la/\mu}(q,t)=Q_{\la/\mu}(1) \\[1mm]
\sk_{\la/\mu}(q,t)&:=Q_{\la/\mu}\Big(\frac{1-q}{1-t}\Big) \\
\widehat{\sk}_{\la/\mu}(q,t)&:=
Q_{\la/\mu}\Big(\frac{1-q/t}{1-t}\Big) \\[1mm]
\ks_{\la/\mu}(q,t)&:=Q_{\la/\mu}(-1),
\end{align}
\end{subequations}
where $\psi'_{\la/\mu}(q,t)$ and $\varphi_{\la/\mu}(q,t)$ is notation
used by Macdonald, and where the $-1$ in $Q_{\la/\mu}(-1)$ is a plethystic 
$-1$, i.e., applied to the power sum $p_r$ of positive index $r$ it gives 
the number $-1$.
The Pieri coefficients $\vs_{\la/\mu}(q,t)$ and $\hs_{\la/\mu}(q,t)$ have
nice factorised forms generalising \eqref{vs} and \eqref{hs}, see 
\cite[pp. 336--342]{Macdonald}. So does $\widehat{\sk}_{\la/\mu}(q,t)$
\cite[p. 173, Remark 2]{Rains06}, \cite[Proposition 3.2]{W05}:
\begin{equation}\label{RW}
\widehat{\sk}_{\la/\mu}(q,t)=\begin{cases}\displaystyle
t^{n(\la)-n(\mu)}
\prod_{i,j=1}^{l(\la)} 
\frac{(qt^{j-i-1};q)_{\la_i-\mu_j} (qt^{j-i};q)_{\mu_i-\mu_j}}
{(qt^{j-i-1};q)_{\mu_i-\mu_j}(qt^{j-i};q)_{\la_i-\mu_j}}
& \text{for $\mu\subseteq\la$}, \\
0 & \text{otherwise},
\end{cases}
\end{equation}
where $(a;q)_k:=(a;q)_{\infty}/(aq^k;q)_{\infty}$ for all $k\in\mathbb{Z}$.
We leave it to the reader to verify that the above right-hand side
for $q=0$ reduces to the right-hand side of \eqref{sk}.
The remaining two Pieri coefficients do not factor into binomials. 
For example
\begin{align*}
\sk_{(2,1)/(1,0)}(q,t)&=\frac{1-q-q^2+t+qt-q^2t}{1-q^2 t} \\
\ks_{(2,1)/(1,0)}(q,t)&=\frac{(1-t)(1+q-t+qt-t^2-qt^2)}{(1-q)(1-q^2 t)}.
\end{align*}
Of course, $\sk_{\la/\mu}(0,t)=\sk_{\la/\mu}(t)$ so it does
factorise in the classical limit. This is however not the case for 
$\ks_{\la/\mu}(0,t)$, and
\[
\ks_{(2,1)/(1,0)}(0,t)=(1-t)(1-t-t^2).
\]

\medskip

Let $g_r=g_r(X;q,t)=Q_{(r)}(X;q,t)$, so that $g_r(X;0,t)=q_r(X;t)$.
Then the following $q$-analogue of Theorem~\ref{thmKL} holds.
\begin{theorem}
For partitions $\nu\subseteq\mu$,
\begin{subequations}
\begin{align}\label{PM1}
P_{\mu/\nu} e_r&=\sum_{\la,\eta}(-1)^{\abs{\nu-\eta}}
\vs_{\la/\mu}(q,t) \sk_{\nu/\eta}(q,t) P_{\la/\eta} \\
P_{\mu/\nu} h_r&=\sum_{\la,\eta}(-1)^{\abs{\nu-\eta}}
\sk_{\la/\mu}(q,t) \vs_{\nu/\eta}(q,t) P_{\la/\eta}
\label{PM2} \\
P_{\mu/\nu}g_r&= \sum_{\la,\eta} 
\hs_{\la/\mu}(q,t) \ks_{\nu/\eta}(q,t) P_{\la/\eta}
\label{PM3} \\
&= \sum_{\la,\eta,\omega} 
(-1)^{\abs{\nu-\omega}} t^{\abs{\omega-\eta}}  
\hs_{\la/\mu}(q,t) \vs_{\nu/\omega}(q,t)\,\widehat{\sk}_{\omega/\eta}(q,t)
P_{\la/\eta},
\label{PM4}
\end{align}
\end{subequations}
where each of the multiple sums is subject to the restriction
$\abs{\la}+\abs{\eta}=\abs{\mu}+\abs{\nu}+r$.
\end{theorem}

\section{The $q$-binomial theorem for Macdonald polynomials}

In \cite[Equation (2.11)]{LW11} Lascoux and the author proved the following
$q$-binomial theorem for Macdonald polynomials:
\begin{equation}\label{qbt}
\sum_{\la} Q_{\la/\nu}\Big(\frac{a-b}{1-t}\Big) 
P_{\la/\mu}(X)
=\bigg(\prod_{x\in X} \frac{(bx;q)_{\infty}}{(ax;q)_{\infty}}\bigg)
\sum_{\la} Q_{\mu/\la}\Big(\frac{a-b}{1-t}\Big) 
P_{\nu/\la}(X).
\end{equation}
For $\mu=\nu=0$ and $(a,b)\mapsto (1,a)$ this is the well-known
Kaneko--Macdonald $q$-binomial theorem \cite{Kaneko96,Macdonald}
\begin{equation}\label{KM}
\sum_{\la} t^{n(\la)} \frac{(a)_{\la}}{c'_{\la}}\,
P_{\la}(X)
=\prod_{x\in X} \frac{(ax;q)_{\infty}}{(x;q)_{\infty}},
\end{equation}
where we have used that \cite[p. 338, (6.17)]{Macdonald95}
\[
Q_{\la}\Big(\frac{1-a}{1-t}\Big)=
t^{n(\la)} \frac{(a)_{\la}}{c'_{\la}}.
\]
Here $(a)_{\la}=(a;q,t)_{\la}:=\prod_{i\geq 1} (at^{1-i};q)_{\la_i}$
and $c'_{\la}=c'_{\la}(q,t)$ is the generalised hook polynomial
$c'_{\la}=\prod_{s\in\la}\big(1-q^{a(s)+1}t^{l(s)}\big)$ with
$a(s)$ and $l(s)$ the arm-length and leg-length of the square $s\in\la$.

To show that \eqref{qbt} encodes the skew Pieri formulas 
\eqref{PM1}--\eqref{PM4} we first consider the $\mu=0$ case
\begin{equation}\label{qbt2}
\sum_{\la} Q_{\la/\nu}\Big(\frac{a-b}{1-t}\Big) 
P_{\la}(X)
=P_{\nu}(X)\prod_{x\in X} \frac{(bx;q)_{\infty}}{(ax;q)_{\infty}}.
\end{equation}
If we multiply this by $Q_{\nu/\mu}\big((b-a)/(1-t)\big)$ and
sum over $\nu$ using \eqref{qbt2} with $(\la,\nu,a,b)\mapsto(\nu,\mu,b,a)$
we obtain
\[
\sum_{\la,\nu} 
Q_{\la/\nu}\Big(\frac{a-b}{1-t}\Big) 
Q_{\nu/\mu}\Big(\frac{b-a}{1-t}\Big) 
P_{\la}(X)=P_{\mu}(X).
\]
This implies the orthogonality relation (implicit in \cite{Rains05}
and given in its more general nonsymmetric form in 
\cite[Equation (6.5)]{LRW09}) 
\begin{equation}\label{ortho}
\sum_{\nu} 
Q_{\la/\nu}\Big(\frac{a-b}{1-t}\Big) 
Q_{\nu/\mu}\Big(\frac{b-a}{1-t}\Big) =\delta_{\la\mu}.
\end{equation}
Thanks to \eqref{ortho}, identity \eqref{qbt} is equivalent to 
\[
\sum_{\la,\eta} Q_{\nu/\eta}\Big(\frac{a-b}{1-t}\Big) 
Q_{\la/\mu}\Big(\frac{b-a}{1-t}\Big)P_{\la/\eta}(X) 
=P_{\mu/\nu}(X)
\prod_{x\in X} \frac{(ax;q)_{\infty}}{(bx;q)_{\infty}}.
\]
There are now three special cases to consider.
First, if $b=aq$ then
\[
P_{\mu/\nu}(X) \prod_{x\in X} (1-ax)=
\sum_{\la,\eta} a^{\abs{\la-\mu}+\abs{\nu-\eta}}
Q_{\la/\mu}\Big(\frac{q-1}{1-t}\Big)
Q_{\nu/\eta}\Big(\frac{1-q}{1-t}\Big) 
P_{\la/\eta}(X).
\]
Equating coefficients of $(-a)^r$ and using definition \eqref{Pcoeff}
yields \eqref{PM1}.
Next, if $a=bq$
\[
P_{\mu/\nu}(X) \prod_{x\in X} \frac{1}{1-bx} =
\sum_{\la,\eta} b^{\abs{\la-\mu}+\abs{\nu-\eta}}
Q_{\la/\mu}\Big(\frac{1-q}{1-t}\Big)
Q_{\nu/\eta}\Big(\frac{q-1}{1-t}\Big) 
P_{\la/\eta}(X).
\]
Equating coefficients of $b^r$ and again using \eqref{Pcoeff}
yields \eqref{PM2}.
Finally, if $a=bt$
\[
P_{\mu/\nu}(X)
\prod_{x\in X} \frac{(btx;q)_{\infty}}{(bx;q)_{\infty}}=
\sum_{\la,\eta} b^{\abs{\la-\mu}+\abs{\nu-\eta}}
Q_{\la/\mu}(1) Q_{\nu/\eta}(-1) P_{\la/\eta}(X), 
\]
Equating coefficients of $b^r$ and using \eqref{Pcoeff} gives
\eqref{PM3}. To show that \eqref{PM3} and \eqref{PM4} are
equivalent, we recall Rains' $q$-Pfaff--Saalsch\"utz summation 
for Macdonald polynomials \cite[Corollary 4.9]{Rains05}:
\begin{equation}\label{Eric}
\sum_{\nu} \frac{(a)_{\nu}}{(c)_{\nu}}\,
Q_{\la/\nu}\Big(\frac{a-b}{1-t}\Big)
Q_{\nu/\mu}\Big(\frac{b-c}{1-t}\Big) =
\frac{(a)_{\mu}(b)_{\la}}{(b)_{\mu}(c)_{\la}}\,
Q_{\la/\mu}\Big(\frac{a-c}{1-t}\Big),
\end{equation}
which for $c=a$ is \eqref{ortho}.
Setting $b=a/q$ and $c=a/t$ and using \eqref{Pcoeff} yields
\[
\ks_{\la/\mu}(q,t) =
(t/q)^{\abs{\la-\mu}}  
\frac{(a/q)_{\mu}(a/t)_{\la}}{(a)_{\mu}(a/q)_{\la}}
\sum_{\nu} (-1)^{\abs{\la-\nu}} \frac{(a)_{\nu}}{(a/t)_{\nu}}\,
\vs_{\la/\nu}(q,t)\,\widehat{\sk}_{\nu/\mu}(q,t).
\]
Taking the $a\to\infty$ limit this further simplifies to
\[
\ks_{\la/\mu}(q,t)=\sum_{\nu} (-1)^{\abs{\la-\nu}} t^{\abs{\nu-\mu}}  
\vs_{\la/\nu}(q,t)\,\widehat{\sk}_{\nu/\mu}(q,t),
\]
which proves the equality between \eqref{PM3} and \eqref{PM4}.

\medskip

To conclude let us mention that all other identities of \cite{KL12}
admit simple $q$-analogues. For example,
if we take \eqref{Eric} and specialise $b=a/q$ and $c=at$ then
\[
\sum_{\mu} \frac{(a)_{\mu}}{(at)_{\mu}}\,
(-1)^{\abs{\la-\mu}} \vs_{\la/\mu}(q,t)
Q_{\mu/\nu}\Big(\frac{1-qt}{1-t}\Big)  =
\frac{(a)_{\nu}(a/q)_{\la}}{(a/q)_{\nu}(at)_{\la}}\,
q^{\abs{\la-\nu}} \hs_{\la/\nu}(q,t).
\]
Letting $a\to\infty$ this reduces to
\[
\sum_{\mu} (-t)^{\abs{\la-\mu}} \vs_{\la/\mu}(q,t)
Q_{\mu/\nu}\Big(\frac{1-qt}{1-t}\Big)=\hs_{\la/\nu}(q,t).
\]
For $q=0$ this is \cite[Lemma 5]{KL12}
\[
\sum_{\mu} (-t)^{\abs{\la-\mu}} \vs_{\la/\mu}(t) \sk_{\mu/\nu}(t)=
\hs_{\la/\nu}(t).
\]
Similarly, according to \cite[Equation (6.23)]{LRW09}
\begin{equation}\label{som}
\sum_{\nu} t^{n(\nu)} \frac{(a)_{\nu}}{c'_{\nu}}\, f_{\mu\nu}^{\la}(q,t)=
Q_{\la/\mu}\Big(\frac{1-a}{1-t}\Big).
\end{equation}
For $a=q=0$ this is \cite[Corollary 6]{KL12}
\[
\sum_{\nu} t^{n(\nu)} f_{\mu\nu}^{\la}(t)=\sk_{\la/\mu}(t).
\]

Finally, to obtain a $q$-analogue of \cite[Theorem 7]{KL12} we have to work 
a little harder. First note that
\begin{align} 
P_{\nu}(X) e_m(X) \sum_{r=0}^{\infty} h_r(X)&=
\sum_{\mu} \sk_{\mu/\nu}(q,t) P_{\mu}(X) e_m(X) \notag \\
&=\sum_{\mu} \sum_{\substack{\la \\ \abs{\la-\mu}=m}}
\vs_{\la/\mu}(q,t) \sk_{\mu/\nu}(q,t) P_{\la}(X).
\label{vorm1}
\end{align}
To compute this in a different way, observe that if we set
$a=q$ in \eqref{KM} then
\[
\sum_{\la} t^{n(\la)} \frac{(q)_{\la}}{c'_{\la}}\,
P_{\la}(X)
=\prod_{x\in X} \frac{1}{1-x}=\sum_{r=0}^{\infty} h_r(X).
\]
Using this as well as $e_m=P_{(1^m)}$ we get
\begin{align*}
P_{\nu}(X) e_m(X) \sum_{r=0}^{\infty} h_r(X)&=
\sum_{\eta} t^{n(\eta)} \frac{(q)_{\eta}}{c'_{\eta}}\,
P_{\nu}(X) P_{\eta}(X) P_{(1^m)}(X).
\end{align*}
By a double use of 
$P_{\mu} P_{\nu} = f_{\mu\nu}^{\la} P_{\la}$ 
this leads to
\begin{align}
P_{\nu}(X) e_m(X) \sum_{r=0}^{\infty} h_r(X)
&=\sum_{\eta} t^{n(\eta)} \frac{(q)_{\eta}}{c'_{\eta}}\,
P_{\nu}(X) P_{\eta}(X) P_{(1^m)}(X) \notag \\
&=\sum_{\mu,\eta} t^{n(\eta)} \frac{(q)_{\eta}}{c'_{\eta}}\,
f_{\eta,(1^m)}^{\mu}(q,t) P_{\mu}(X) P_{\nu}(X) \notag \\
&=\sum_{\la,\mu,\eta} t^{n(\eta)} \frac{(q)_{\eta}}{c'_{\eta}}\,
f_{\eta,(1^m)}^{\mu}(q,t) f_{\mu\nu}^{\la}(q,t) P_{\la}(X) \notag \\
&=\sum_{\la,\mu} \sk_{\mu/(1^m)}(q,t) f_{\mu\nu}^{\la}(q,t) P_{\la}(X),
\label{vorm2}
\end{align}
where the final equality follows from the $a=q$ case of \eqref{som}.
Equating coefficient of $P_{\la}(X)$ in \eqref{vorm1} and \eqref{vorm2}
yields
\[
\sum_{\substack{\mu \\ \abs{\la-\mu}=m}}
\vs_{\la/\mu}(q,t) \sk_{\mu/\nu}(q,t) =
\sum_{\mu} \sk_{\mu/(1^m)}(q,t) f_{\mu\nu}^{\la}(q,t) .
\]
By \eqref{sk},
\[
\sk_{\la/(1^m)}(0,t)=\sk_{\la/(1^m)}(t)=
t^{n(\la/(1^m))} \qbin{\la'_1}{m}_t=
t^{n(\la)-\binom{m}{2}} \qbin{\la'_1}{m}_{t^{-1}},
\]
so that for $q=0$ we obtain \cite[Theorem 7]{KL12}
\[
\sum_{\substack{\mu \\ \abs{\la-\mu}=m}}
\vs_{\la/\mu}(t) \sk_{\mu/\nu}(t) =
\sum_{\mu} t^{n(\la)-\binom{m}{2}} 
f_{\mu\nu}^{\la}(t) \qbin{\la'_1}{m}_{t^{-1}}.
\]

\subsection*{Acknowledgements}
I thank Matja\v{z} Konvalinka and Aaron Lauve for helpful discussions.

\bibliographystyle{amsplain}

\end{document}